\documentclass{amsart}

\usepackage{graphicx,color}
\usepackage[initials]{amsrefs}
\usepackage[all]{xy}

\makeatletter
\renewcommand\MR[1]{%
    \relax\ifhmode\unskip\spacefactor3000 \space\fi
    MR\thinspace \MRno#1 \@nil
}
\def\MRno MR#1 #2\@nil{\@stripzeroes#1}
\def\@stripzeroes{%
  \@ifnextchar0{\expandafter\@stripzeroes\@gobble}{}}
\makeatother

\def\pres<#1;#2>{\left<\,#1 \mathrel; #2 \,\right>}
\newcommand\abs[1]{\lvert#1\rvert}
\newcommand\iso{\cong}
\newcommand\Z{\ensuremath{\mathbb{Z}}}
\newcommand\into{\hookrightarrow}
\newcommand\Q{\ensuremath{\mathbb{Q}}}
\newcommand\B{\ensuremath{\mathcal{B}}}
\newcommand\F{\ensuremath{\mathcal{F}}}
\newcommand\M{\ensuremath{\mathcal{M}}}
\renewcommand\S{\ensuremath{\mathcal{S}}}
\newcommand\ab{\text{ab}}

\theoremstyle{plain}
\newtheorem{theorem}{Theorem}[section]

\theoremstyle{remark}
\newtheorem{remark}[theorem]{Remark}

\numberwithin{equation}{section}

\begin{document}

\title{The Homology of Richard Thompson's Group $F$}

\author{Kenneth S. Brown}

\address{Department of Mathematics\\ 
Cornell University\\ 
Ithaca, NY 14853}
\dedicatory{Dedicated to Ross Geoghegan in honor of his
  60th birthday.}
\email{kbrown@math.cornell.edu}

\subjclass[2000]{Primary 20J06; Secondary 55P20}

\keywords{Cohomology of groups, Thompson's group, piecewise linear
homeomorphisms, Cantor algebras}

\maketitle

\section*{Introduction}
\label{sec:introduction}

Let $F$ be Richard Thompson's group, which can be defined by the
presentation
\[
F = \pres<x_0,x_1,x_2,\dots; x_n^{x_i}=x_{n+1} \mbox{ for } i<n>\;.
\]
Here $x^y:= y^{-1}xy$.  One can also describe $F$ in a variety of
other ways, some of which are reviewed briefly in
Section~\ref{sec:background-f}.  In the early 1980s Ross Geoghegan and
I studied the homological properties of~$F$;
see~\cites{brown85:_cohom,brown84:_fp,brown83:_fp_hnn}.  We showed,
among other things, that the integral homology~$H_*(F)$ is free
abelian of rank~2 in every positive dimension.  It turns out that the
homology admits a natural ring structure, which I calculated a few
years later.  The answer is that $H_*(F)$ is an associative ring
(without identity) generated by an element~$\varepsilon$ of degree~0 and
elements $\alpha,\beta$ of degree~1, subject to the relations
\begin{gather*}
\varepsilon^2=\varepsilon\\
\varepsilon\alpha=\beta \varepsilon=0\\
\alpha \varepsilon =\alpha\;, \quad \varepsilon\beta=\beta\;.
\end{gather*}
It follows that $\alpha^2=\beta^2=0$ and that the alternating products
$\alpha\beta\alpha\cdots$ and $\beta\alpha\beta\cdots$ give a basis
for the homology in positive dimensions.

With the aid of this ring structure on the homology, one can calculate
the integral cohomology ring:
\begin{equation}
\label{eq:1}
H^*(F) \cong \bigwedge(a,b) \otimes \Gamma(u)\;,
\end{equation}
where $\bigwedge(a,b)$ is an exterior algebra on two generators~$a,b$
of degree~1, and $\Gamma(u)$ is a divided polynomial ring on one
generator~$u$ of degree~2.

I never published these results because they were subsumed by the work
of Melanie Stein~\cite{stein92:_group}, who proved analogous results
for a much larger class of groups.  Since $F$ remains of great current
interest, and since readers may find it inconvenient to work through
Stein's paper and specialize everything to the case of~$F$, I give
here my original proofs for that case.  In particular, I explain the
cohomology calculation~\eqref{eq:1}, which is not stated explicitly
in~\cite{stein92:_group}.

The impetus for publishing these 15-year-old results at the present
time comes from a question recently raised by Geoghegan [private
communication]:  Is $F$ a ``K\"ahler group'', i.e., the fundamental
group of a K\"ahler manifold?  Now one of the necessary conditions for
a group to be a K\"ahler group is that the cup product
\mbox{$H^1\otimes H^1
\to H^2$} be nontrivial.  So the calculation~\eqref{eq:1} is consistent
with an affirmative answer to Geoghegan's question.  In fact, $F$
satisfies all of the necessary conditions for being a K\"ahler group
that I am aware of.

The remainder of the paper is organized as follows.
Section~\ref{sec:background-f} is a brief review of~$F$ and some of
its properties.  I introduce in Section~\ref{sec:product-f} a
homomorphism \mbox{$\mu\colon F\times F\to F$}, which induces the product
structure on~$H_*(F)$ alluded to above.  The product~$\mu$ is not
strictly associative, but it is associative up to conjugacy, so the
homology becomes an associative ring.  In
Section~\ref{sec:cubical-kf-1} I describe a cubical
$K(F,1)$-complex~$X$ due to Stein.  It is a variant of a $K(F,1)$
constructed in~\cite{brown84:_fp}, but it has the advantage that it
comes equipped with a strictly associative product $X\times X\to X$,
which induces the product~$\mu$ on~$\pi_1$.  [Note: One has to be
careful about basepoints in order to make sense of the last statement,
since the product is not basepoint-preserving.]  The cellular chain
complex~$C_*(X)$ then becomes a differential graded ring.  Its
structure is so simple that one can compute its homology directly;
this is done in Section~\ref{sec:calc-homol-ring}.  Finally, I explain
in Section~\ref{sec:calc-cup-prod} how the homology calculation yields
the cohomology result stated in~\eqref{eq:1}.

This paper is an expanded version of a talk I gave at an AMS special
session in Nashville on October~17, 2004, in honor of Ross Geoghegan's
60th birthday.  Ross is the person who first told me about Thompson's
group and made me realize what a fascinating object it is.  It is a
great pleasure to dedicate this paper to him.

\section{Background on $F$}
\label{sec:background-f}

There are many references for the basic facts about~$F$, including
\cites{belk04:_f,belk:_fores_f,brown87:_finit,brown84:_fp,cannon96:_introd_richar},
all of which have further references in their bibliographies.
We review in this section only a few basic facts that we will need
later.

\subsection{Dyadic PL-homeomorphisms}
\label{sec:dyad-pl-home}

\begin{figure}
\centering
\includegraphics{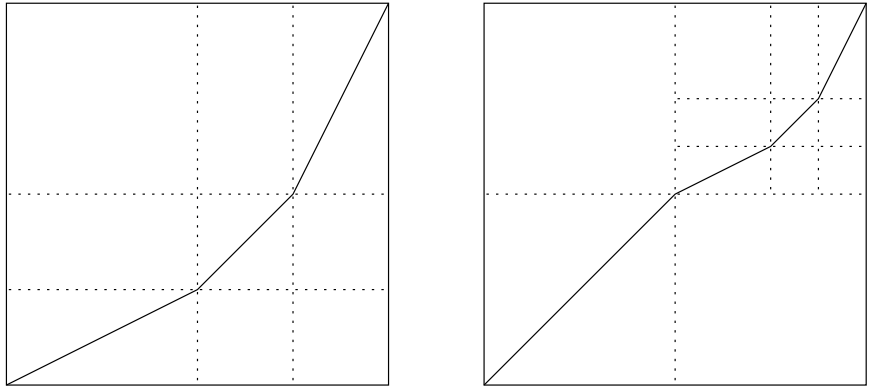}
\caption{The generators $x_0$ and $x_1$}
\label{fig:homeo}
\end{figure}

Let $I$ and $J$ be intervals of real numbers.  A homeomorphism
$f\colon I\to J$ will be called \emph{dyadic} if it satisfies the
following conditions:
\begin{itemize}
\item
$f$ is piecewise linear with only finitely many breakpoints.
\item
Every breakpoint has dyadic rational coordinates.
\item
Every slope is an integral power of~2.
\end{itemize}
A basic fact about $F$, essentially known to Thompson, is that $F$ can
be identified with the group of dyadic homeomorphisms of the unit
interval~$[0,1]$, with the group law being composition:  $(fg)(t) =
f(g(t))$ for $f,g\in F$ and $t\in [0,1]$.  [Warning:  Some authors,
including Thompson himself, have used the opposite convention for
composition.]  Figure~\ref{fig:homeo} shows the graphs of the first
two generators, $x_0$ and~$x_1$.
Note that $x_1$ is the identity on~$[0,1/2]$ and is a rescaled copy
of~$x_0$ on~$[1/2,1]$.  The remaining generators $x_2,x_3,\dots$ are
obtained by repeating this process.  More precisely, there is a
homomorphism $\phi\colon F\to F$ such that $\phi(f)$ is the identity
on~$[0,1/2]$ and is a rescaled copy of~$f$ on~$[1/2,1]$, and the
assertion is that $\phi(x_n)=x_{n+1}$ for all $n\geq 0$.

All dyadic PL-homeomorphisms of $[0,1]$, and hence all elements
of~$F$, can be obtained from \emph{dyadic subdivisions}, as
illustrated by the dotted lines in Figure~\ref{fig:homeo}.  One
subdivides the domain and range into an equal number of parts by
repeated insertion of midpoints, and one maps the subintervals of the
domain subdivision linearly to the subintervals of the range
subdivision.

\begin{figure}[htb]
\centering
\includegraphics{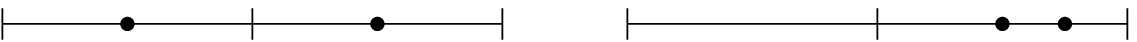}
\caption{Two dyadic subdivisons of $[0,2]$}
\label{fig:subdivisions}
\end{figure}

\begin{figure}[htb]
\centering
\includegraphics{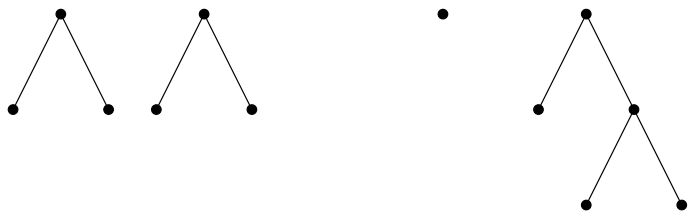}
\caption{Forests represent dyadic subdivisions}
\label{fig:forests}
\end{figure}

Finally, we recall for future reference that there is a fairly obvious
way to represent dyadic subdivisions by rooted binary trees, as
explained in
\cites{belk04:_f,belk:_fores_f,brown87:_finit,stein92:_group}.  More
generally, we will have occasion to talk about dyadic subdivisions of
$[0,n]$, where $n$ is an integer~$\geq 1$, in which we start with the
subdivision into unit intervals $[i-1,i]$ ($i=1,\dots,n$) and then
dyadically subdivides these, as in 
Figure~\ref{fig:subdivisions}.
Such a subdivision is represented by a binary forest consisting of $n$
rooted binary trees (in a definite order).  Thus the roots of the
forest correspond to the intervals $[i-1,i]$, and the leaves
correspond to the parts of the sudivision.  For example, the two
subdivisions of $[0,2]$ in Figure~\ref{fig:subdivisions} are
represented by the forests in Figure~\ref{fig:forests}.

\subsection{Conjugacy idempotents}
\label{sec:conj-idemp}

Dydak~\cites{dydak77:_fanr,dydak77} and, independently, Freyd and
Heller~\cite{freyd93:_split} rediscovered Thompson's group in
connection with a problem in homotopy theory.  Briefly, they were
interested in ``free homotopy idempotents'', i.e., maps from a space
to itself that are idempotent up to homotopy; the maps are required to
be basepoint preserving, but the homotopy is not.  Passing to the
fundamental group, one is led to study group homomorphisms that are
idempotent up to conjugacy, and it turns out that $F$ carries the
universal example of such a homomorphism.  Namely, the defining
relations for~$F$ show that the shift homomorphism~$\phi$ discussed
above satisfies
\[
\phi^2(f) = \phi(f)^{x_0}\;
\]
for all $f\in F$, so $\phi$ is idempotent up to conjugacy.  To see
universality, suppose we have a group~$G$ and an endomorphism
$\psi\colon G\to G$ such that $\psi^2 = \psi^{y_0}$ for some $y_0 \in
G$.  Then $G$ contains elements $y_n:=\psi^n(y_0)$ for $n\geq 1$.  The
equation $\psi^2=\psi^{y_0}$ implies that $y_j^{y_i}= y_{j+1}$ for
$j>i=0$, and we can repeatedly apply~$\psi$ to see that this remains
valid for $j>i\geq 0$.  Thus there is a homomorphism $F\to G$ taking
$x_n$ to~$y_n$ for all~$n$.

\subsection{Algebra automorphisms}
\label{sec:algebra-autom}

Consider the algebraic system consisting of a set~$A$ together with a
bijection $A\times A\to A$.  Thus there is a product on~$A$, and every
$a\in A$ factors uniquely as a product of two other elements.
Algebras of this type have appeared in several places; see
\cites{smirnov71:_cantor,swierczkowski61,jonsson61,clark69:_variet,higman74:_finit}.
Following Smirnov~\cite{smirnov71:_cantor}, we call~$A$ a \emph{Cantor
algebra}.  Galvin and Thompson [unpublished] showed that $F$ is
isomorphic to the group of ``order-preserving'' automorphisms of the
free Cantor algebra on one generator.  This point of view was
exploited in~\cite{brown87:_finit}.  In order to explain what
``order-preserving'' means, we need to recall some facts about bases
of free Cantor algebras.  The definitions will be adapted to our
present needs, in which bases are always ordered.

All Cantor algebras considered in this paper will come equipped with a
``standard basis'' $a_1,\dots,a_n$, in which the order of the basis
elements is important.  A \emph{simple expansion} of an ordered basis
consists of factoring one of the basis elements~$a$ as~$a_0a_1$ and
replacing $a$ by the two elements $a_0,a_1$ (in that order).  A
general expansion consists of doing finitely many simple expansions.
It is an easy fact, proved in several of the references cited above,
that an expansion of a basis is again a basis.  For example, if $A$ is
the free Cantor algebra on one generator~$a$, then we can factor~$a$
as a triple product in two different ways:
\[
a=a_0a_1=a_0(a_{10}a_{11})\;, \qquad a=a_0a_1=(a_{00}a_{01})a_1\;.
\]
This yields a basis $a_0,a_1$ of size~2 and two bases
$a_0,a_{10},a_{11}$ and $a_{00},a_{01},a_1$ of size~3.  The opposite
of expansion is called \emph{contraction}.  A \emph{simple
contraction} consists of replacing two consecutive basis elements by
their product, and a general contraction consists of finitely many
simple contractions.  Finally, an \emph{ordered basis} of our Cantor
algebra is one that can be obtained from the standard basis by doing
finitely many expansions or contractions.

\begin{remark}
\label{rem:2}
Expanding a basis is analogous to subdividing an interval.  In
particular, if we start with an ordered basis having $n$ elements,
then there is a fairly obvious way to represent a $k$-fold expansion
of it by a binary forest with $n$ roots and $n+k$ leaves.  For
example, given an ordered basis $a,b$, we have a 2-fold expansion
$a_0,a_1,b_0,b_1$ and a second 2-fold expansion $a,b_0,b_{10},b_{11}$.
These are represented by the two forests in Figure~\ref{fig:forests}
above.
\end{remark}

Given two Cantor algebras of the type we are considering (free with a
given linearly ordered basis), we define an \emph{order-preserving}
isomorphism between them to be one that takes an ordered basis of the
domain to an ordered basis of the range (preserving the ordering on
the bases).  For example, if $A$ is free on~$a$ as above, then there
is an automorphism of~$A$ taking $a_0,a_{10},a_{11}$ to
$a_{00},a_{01},a_1$.  This corresponds to~$x_0$.  Using the analogy
between subdivision and expansion, the reader should be able to look
at Figure~\ref{fig:homeo} and guess which automorphism of~$A$
corresponds to~$x_1$.

The ordered bases of the free Cantor algebra~$A$ on one generator form
a poset~\B, in which $B \leq C$ if $C$ is an expansion of~$B$.  This
poset is a directed set; it played an important role
in~\cite{brown87:_finit} and will be referred to again later.

\begin{remark}
\label{rem:1}
For readers who prefer to avoid Cantor algebras, here is an alternate
description of the poset \B\ of ordered bases.  Consider dyadic
PL-\hspace{0pt}homeo\-morphisms $f\colon [0,n] \to [0,1]$ where $n$ is
an integer~$\geq 1$.  These play the role of bases.  [To see why, note
that giving an ordered basis of~$A$ of size~$n$ is equivalent to
giving an order-preserving isomorphism $A_n \to A$, where $A_n$ is
free on $a_1,\dots,a_n$.]  An \emph{expansion} of~$f$ is a
homeomorphism $[0,n+k]\to [0,1]$ of the form $f\circ s$, where
$s\colon [0,n+k] \to [0,n]$ is a \emph{subdivision map}, obtained as
follows:  Perform a dyadic subdivision of $[0,n]$ into $n+k$ parts,
and let $s$ map the unit intervals $[j-1,j]$ of~$[0,n+k]$ to the $n+k$
subintervals of~$[0,n]$.  See Figure~\ref{fig:subdivision_map} for an
example with $n=1$ and $k=2$.
\begin{figure}[htb]
\centering
\includegraphics{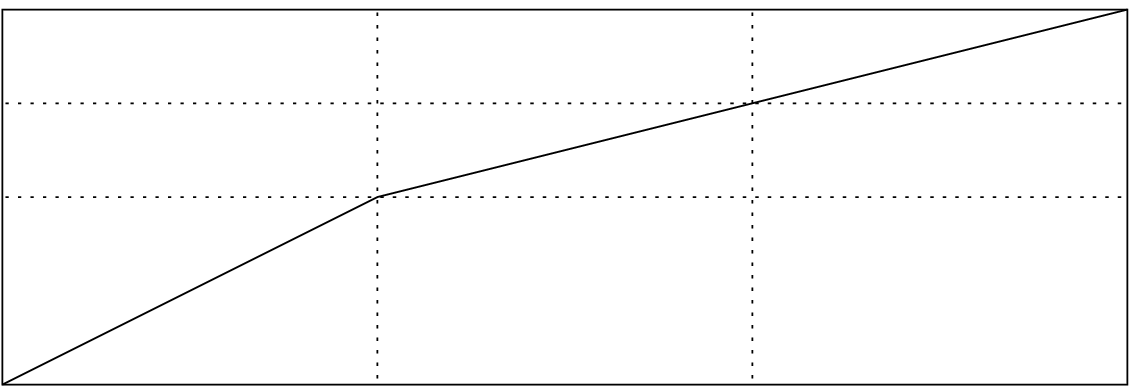}
\caption{A subdivision map $[0,3]\to[0,1]$}
\label{fig:subdivision_map}
\end{figure}
\end{remark}

\subsection{Finiteness properties of $F$}
\label{sec:finit-prop-f}

It is obvious that $F$ is generated by the first two generators
$x_0,x_1$, since the remaining~$x_n$ are obtained by repeated
conjugation.  Less obviously, two relations suffice:
\begin{align*}
x_1^{x_0x_0} &= x_1^{x_0x_1}\\[\jot]
x_1^{x_0x_0x_0} &= x_1^{x_0x_0x_1}\;.
\end{align*}
This pattern was extended in Brown--Geoghegan~\cite{brown84:_fp},
where it was shown that there is an Eilenberg--MacLane complex~$Y$ of
type $K(F,1)$ with exactly two cells in each positive dimension.  Our
method was motivated by the connection between $F$ and the theory of
free homotopy idempotents.  Namely, we constructed the universal
example of a space~$X$ with a free homotopy idempotent, and we showed
that $X$ was a $K(F,1)$-complex.  Now $X$ has infinitely many cells in
each positive dimension, but we were able to show, by imitating the
proof that $F$ requires only two generators and two relations, that
there were only two ``homotopically essential'' cells in each
dimension.  The desired complex~$Y$ was then obtained as a quotient
complex of~$X$.

Incidentally, this is where the homology calculation cited in the
introduction came from:  We showed that the cellular chain
complex~$C_*(Y)$ has trivial boundary operator, so that $H_n(F)\iso
\Z^2$ for each $n\geq 1$.

\subsection{The abelianization of $F$}
\label{sec:abelianization-f}

Let $F'$ be the commutator subgroup of~$F$ and let $F_{\ab}$ be
the abelianization~$F/F'$; it is isomorphic to~$\Z\times\Z$.  From the point
of view of dyadic PL-homeomorphisms, the abelianization map $F\to
\Z\times\Z$ is given by $f\mapsto \bigl(\log_2 f'(0), \log_2
f'(1)\bigr)$.  [Recall that the slopes of elements of~$F$ are in the
infinite cyclic multiplicative group generated by~2.]  Thus $F'$
consists of the elements of~$F$ that are the identity near the
endpoints~$0,1$.  The group~$F'$ is simple but infinitely generated;
it is the union of an increasing sequence of isomorphic copies of~$F$.

\section{A product on $F$}
\label{sec:product-f}

\subsection{Definition}
\label{sec:definition}

\begin{figure}[htb]
\centering
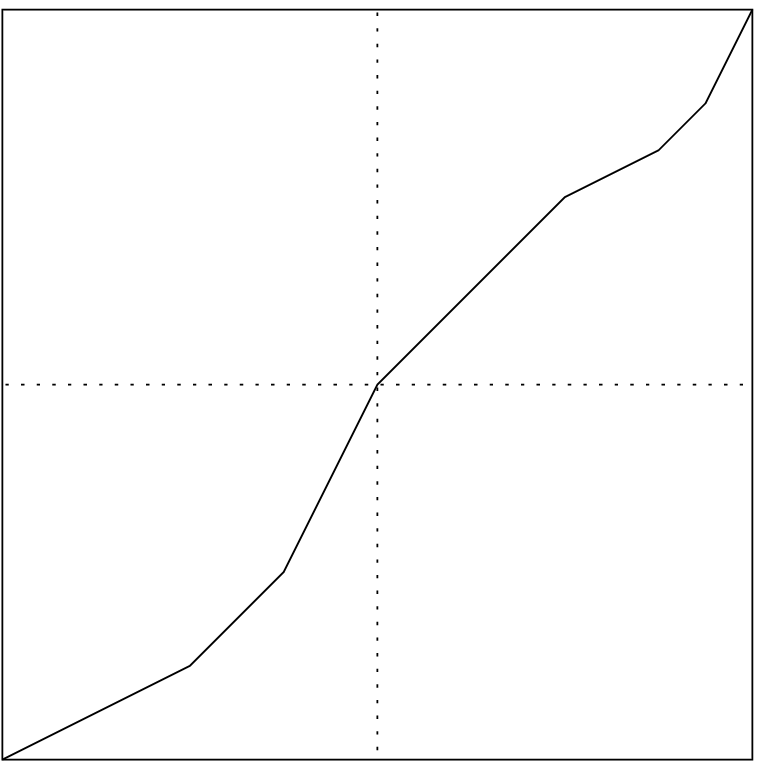
\caption{$x_0*x_1$}
\label{fig:product}
\end{figure}

There is a homomorphism $\mu\colon F\times F\to F$, denoted
$(f,g)\mapsto f*g$, which is defined as follows:  If we interpret $F$
as the group of dyadic homeomorphisms of~$[0,1]$, then $f*g$ is given
by
\[
(f*g)(t) = 
\begin{cases}
f(2t)/2 & 0\leq t \leq 1/2\\[\jot]
(g(2t-1)+1)/2 & 1/2 \leq t \leq 1.
\end{cases}
\]
Less formally, $f*g$ is a rescaled copy of $f$ on $[0,1/2]$ and a
rescaled copy of~$g$ on~$[1/2,1]$; see Figure~\ref{fig:product} for an
example.
Alternatively, if we interpret $F$ as the group of order-preserving
isomorphisms of the free Cantor algebra~$A$ on one generator~$a$, then
we can first form the automorphism $f\amalg g$ of $A\amalg A$, where
$A\amalg A$ is the categorical sum of two copies of~$A$ and hence is
free on two generators, and then we can transport this to an
automorphism of~$A$ via a suitable isomorphism $A\iso A\amalg A$.

One can check that $\mu$ is associative up to conjugacy:
\begin{equation}
\label{eq:2}
f*(g*h) = ((f*g)*h)^{x_0}\;.
\end{equation}
There is, however, no identity.  In particular, the identity element~1
for the group law on~$F$ is not an identity for~$\mu$; indeed, one has
\begin{equation}
\label{eq:3}
1*f=\phi(f)\;,
\end{equation}
where $\phi$ is the conjugacy idempotent discussed above.  Similarly,
$\mu$-multiplying on the right by~1 is a (different) conjugacy
idempotent on~$F$.

\subsection{The induced product on homology}
\label{sec:induc-prod-homol}

Exactly as in the homology theory of abelian groups
\cite{brown94:_cohom}*{V.5}, the product~$\mu$ induces a product
\[
H_*(F) \otimes H_*(F)\to H_*(F)\;,
\]
making $H_*(F)$ an associative ring, where associativity comes
from~\eqref{eq:2} and the fact that inner automorphisms act trivially
on homology \cite{brown94:_cohom}*{II.6.2}.  This product has no
identity; the canonical generator~$\varepsilon$ of $H_0(F)=\Z$ is
idempotent, and left multiplication by it is the endomorphism $\phi_*$
of~$H_*(F)$ by~\eqref{eq:3}.  As we will see, this idempotent
endomorphism has rank~1 in every positive dimension, as does the
encomorphism given by right multiplication by~$\varepsilon$.

We now wish to calculate the product explicitly.  We will do this in
Section~\ref{sec:calc-homol-ring} after recalling a construction due
to Stein~\cite{stein92:_group}.

\section{Stein's cubical $K(F,1)$}
\label{sec:cubical-kf-1}

\subsection{A groupoid analogue of $F$}
\label{sec:groupoid-analogue-f}

Let \F\ be the category whose objects are the intervals $[0,n]$
($n\geq 1$) and whose morphisms are the dydadic PL-\hspace{0pt}homeo\-morphisms
between them.  Alternatively, we could take the objects to be the free
Cantor algebras~$A_n$ and the morphisms to be the order-preserving
isomorphisms.  Let $\abs{\F}$ denote the geometric realization of~\F\
as defined by Quillen~\cite{quillen73:_higher_k}.  See
also~\cite{quillen78:_homot}, where the theory is reviewed for the
case where the category is a poset.  [A poset can be viewed as a
category, with one morphism $a\to b$ for every relation $a\leq b$.]
Recall that $\abs{\F}$ is a CW-complex in which every open $p$-cell can
be identified with the interior of the standard $p$-simplex.  There is
one such $p$-cell for every $p$-tuple of composable non-identity
morphisms in~\F:
\[
\xymatrix{I_{n_0}\ar[r]^{f_1} & I_{n_1}\ar[r]^{f_2} &
  \cdots\ar[r]^{f_p} & I_{n_p}} \;.
\]
Faces are gotten by deleting objects and, if the object is not the
first or last, composing morphisms.  [It is possible that this will
yield an identity map, hence a ``degenerate'' simplex that must be
collapsed to a lower-dimensional simplex.]  Since \F\ is a connected
groupoid and the group of automorphisms of one object is~$F$, it
follows from Quillen's theory that $\abs{\F}$ is an Eilenberg--MacLane
complex of type~$K(F,1)$.  We will not actually need to make use of
this fact, but it provides helpful motivation.

\subsection{A smaller category}
\label{sec:smaller-category}

The complex~\F\ is much too big to be of any use computationally, so
we pass to a subcategory~\S\ such that $\abs{\S}$ is still a $K(F,1)$.
The objects of~\S\ are the same as those of~\F, but as morphisms we
only use the subdivision maps~$s$ as defined in Remark~\ref{rem:1} in
Section~\ref{sec:algebra-autom}.  Alternatively, if we want to use
Cantor algebras as the objects, we only use isomorphisms $A_{n+k} \to
A_n$ that take the standard basis elements of~$A_{n+k}$ to basis
elements of~$A_n$ obtained by doing a $k$-fold expansion of its
standard basis.  From either point of view, one sees immediately that
the morphisms in~\S\ from the object associated with ~$n+k$ to the
object associated with~$n$ correspond to binary forests with $n$ roots
and $n+k$ leaves.  See Belk \cite{belk04:_f}*{Section~7.2} for
further remarks on~\S.

Although we will not need this fact, one can check that $\abs{\S}$ is
a regular CW-complex in which every closed cell is canonically
homeomorphic to a standard simplex.  It is not a simplicial complex,
however, because a simplex is not determined by its vertices.  But we
will see shortly that its universal cover is a simplicial complex.

To see that $\abs{\S}$ is a $K(F,1)$-complex, we can proceed in two
different ways.  The first method, based on
\cite{quillen73:_higher_k}*{Theorem~A}, is to consider the fibers in
the sense of Quillen of the inclusion $\S\into\F$; it suffices to show
that they are contractible.  A direct check of the definitions shows
that each fiber is a poset and is a directed set, hence it is indeed
contractible.  For example, the fiber over~$I_1$ is isomorphic to the
poset~\B\ discussed in Section~\ref{sec:algebra-autom}, with the order
relation reversed.

A more elementary approach is to directly construct the universal
cover of~$\abs{\S}$ and observe that it is contractible.  In fact, the
universal cover turns out to be the geometric realization~$\abs{\B}$
of the contractible poset~\B.  One can see this by thinking of~\B\
either as the poset of bases or as a poset constructed using dydadic
maps of intervals.  From either point of view there is an obvious
action of $F$ on~\B, and it is straightforward to check that the
induced action on the (contractible) simplicial complex~$\abs{\B}$ is
free and that $\abs{\S}$ is the quotient.  Hence $\abs{\S}$ is indeed
a $K(F,1)$-complex.

\subsection{A cubical complex}
\label{sec:cubical-complex}

Following Stein~\cite{stein92:_group}, we now pass to a further
subcomplex~$X\subset \abs{\S}$, which is still a $K(F,1)$ and in which
the simplices can be naturally grouped into cubes.  It is easiest to
carry this out in the universal cover~$\abs{\B}$ and then pass to the
quotient by the action of~$F$.  We will work with the original
definition of \B\ as the poset of ordered bases of~$A=A_1$.  Readers
who prefer to work with intervals and dyadic subdivisions can
translate everything into that language or can refer
to~\cite{stein92:_group}.

Given an ordered basis $L$, an \emph{elementary expansion} of~$L$ is a
basis~$M$ obtained by replacing zero or more elements $b\in L$ by
their factors $b_0,b_1$.  Such an expansion corresponds to an
``elementary forest'', by which we mean one in which each tree
consists either of the root only or the root with a single pair of
descendants.  For example, the expansion represented by the forest on
the left in Figure~\ref{fig:forests} is elementary, while the one on
the right is not.  We will write $L\preceq M$ if $M$ is an elementary
expansion of~$L$.  [Warning:  This relation is not transitive.]
Recall that a simplex of~$\abs{\B}$ is given by a chain
$L_0<L_1<\cdots<L_p$ of bases.  Call the simplex \emph{elementary} if
$L_0\preceq L_p$.  This implies that $L_i\preceq L_j$ for $0\leq i
\leq j \leq p$.  Hence every face of an elementary simplex is
elementary, and the elementary simplices form an $F$-invariant
subcomplex~$\tilde{X}$ of~$\abs{\B}$.  Passing to the quotient by the
action of~$F$, we obtain the desired subcomplex $X\subset \abs{\S}$.
It has one cell for each chain of bases $L_0<\cdots <L_p$ such that
$L_0$ is the standard basis of~$A_n$ for some $n\geq 1$ and $L_p$ is
an elementary expansion of~$L_0$.

Stein proves that the complex of elementary simplices is contractible,
so that $X$ is again a $K(F,1)$-complex.  Her proof is given in detail
in~\cite{stein92:_group} and is repeated, in a slightly different
context, in~\cite{brown92}, so we will not repeat it again here.

Finally, we will explain, still following Stein, how to give a coarser cell
decomposition of~$X$ by lumping certain simplices into
cubes, one for each elementary expansion $L\preceq M$ in which $L$ is
the standard basis of some some~$A_n$.  Once again, it is easier to
first do this in the universal cover, so we consider $L\preceq M$ with
arbitrary $L\in\B$ and with $M$ a $k$-fold elementary expansion
of~$L$.  Then the interval $[L,M]$ in~$\B$ is isomorphic, as a poset,
to~$\{0,1\}^k$, where $0<1$ and the product is ordered component-wise.
Since the geometric realization of $\{0,1\}$ is canonically
homeomorphic to the unit interval~$[0,1]$, we conclude that the
geometric realization~$\bigl|[L,M]\bigr|$ is a simplicially subdivided
$k$-cube.

The relative interior of this $k$-cube is the union of the open
simplices corresponding to the chains $L_0<\cdots<L_p$ with $L_0=L$
and $L_p=M$, so these interiors partition~$\tilde{X}$.  We therefore
have a decomposition of~$\tilde{X}$ as a regular cell complex in which
all the closed cells are cubes.  Passing to the quotient by~$F$, we
obtain the desired decomposition of~$X$.  The closed cells are cubes
with some identifications on the boundary.  For example, there is a
2-cube in~$X$ for every $1\leq i < j\leq n$, corresponding to the
following diagram in~\S:
\[\def\labelstyle{\textstyle}
\xymatrix{A_{n+2} \ar[r]^{s_i} \ar[d]_{s_{j+1}} \ar@{-->}[dr] &
  A_{n+1} \ar[d]^{s_j}\\
A_{n+1} \ar[r]_{s_i} & A_n}
\]
Here the top $s_i$ is the isomorphism that maps the standard basis
of~$A_{n+2}$ to the simple expansion at position~$i$ of the standard
basis of~$A_{n+1}$; the other maps are defined similarly.  The diagram
reflects the fact that there are two ways to expand the standard basis
of~$A_n$ at positions~$i,j$:  First expand at position~$j$ and then
expand the $i$th basis element of the result, or first expand at
position~$i$ and then expand the $(j+1)$st basis element of the result.

\subsection{A product on $X$}
\label{sec:product-x}

There is a product $\mu\colon \F\times\F\to\F$, analogous to the
product on~$F$ defined in Section~\ref{sec:definition}, except that no
rescaling is required.  From the point of view of dyadic
PL-homeomorphisms, we define~$\mu$ on objects by setting $I_n * I_m =
I_{n+m}$, and we define it on maps by gluing them together in the
obvious way.  Thus if $f$ has domain~$I_n$ and $g$ has domain~$I_m$,
then $f*g$ is $f$ on~$[0,n]$ and a shifted copy of~$g$ on~$[n,n+m]$.
Alternatively, in terms of Cantor algebras, we set $A_n * A_m = A_n
\amalg A_m$, which we identify with~$A_{n+m}$, and we set $f*g=f\amalg
g$.

\begin{figure}[htb]
\centering
\includegraphics{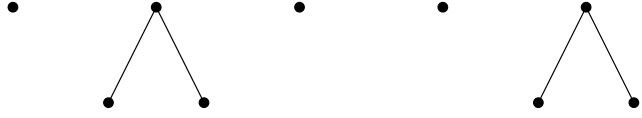}
\caption{The product $vev^2e$}
\label{fig:semigroup}
\end{figure}

Notice that the product $\mu$ is strictly associative; this is the
advantage of working with the groupoid~\F\ instead of the group~$F$.
Thus $\abs{\F}$ becomes a topological semigroup, with $\abs{\S}$
and~$X$ as subsemigroups.  If we use the cubical structure on~$X$, we
find that the product of two cubes is again a cube, so that the set of
cubes is a semigroup.  This semigroup is quite easy to understand:
Recalling that cubes correspond to elementary expansions, which are
represented by elementary forests, the product is given by taking the
disjoint union of the forests.  It follows at once that the semigroup
of cubes is the free semigroup on two generators $v,e$, where $v$ is
the vertex corresponding to the object~$A_1$ and $e$ is the unique
edge in~$X$ from $v$ to~$v^2$.  [Note that the 0-cube~$v$ is
represented by the trivial forest with one root, and the 1-cube~$e$ is
represented by the forest consisting of a root that has two
descendants.]  For example, the 2-cube corresponding to the elementary
expansion shown in Figure~\ref{fig:semigroup} is the product $vev^2e$.
We are using juxtaposition here instead of~$*$ to denote the product
of cubes, since there is no other product under discussion on the set
of cubes.

It is easy to check that the product $X\times X\to X$ induces the
product previously defined on~$F=\pi_1(X)$, if one is careful with
basepoints.  Namely, since $(v,v)\mapsto v^2$, we use the edge~$e$ to
change basepoints so that there is an induced map on~$\pi_1(X) =
\pi_1(X,v)$.

Since $X$ is a topological semigroup, we get a ring structure
on~$H_*(X)= H_*(F)$, which we are now in a position to calculate.

\section{Calculation of the homology ring}
\label{sec:calc-homol-ring}

Let $C=C(X)$ be the cellular chain complex of~$X$ with respect to its
cubical structure.  Then $C$ is a differential graded ring (without
identity).  If we forget the differential, $C$ is the (graded) ring of
noncommuting polynomials without constant term, in two variables~$v,e$
with $v$ in degree~0 and $e$ in degree~1.  The differential is
determined by the formulas
\[
\partial e = v^2 - v\;, \qquad \partial(xy)=\partial x\cdot y +
(-1)^{\deg x} x\cdot \partial y\;.
\]
Note, for the sake of intuition, that $C$ is the universal example of
a differential graded ring having an element $v$ of degree~0 such that
$v^2$ is homologous to~$v$.  This should be compared to the
description of~$F$ in Section~\ref{sec:conj-idemp}.  We will use $C$
to compute $H_*(F)$ as a ring.  What makes the computation feasible is
that, although $C$ is a free abelian group of infinite rank in every
dimension, it is ``small'' in the sense of having only two generators
as a ring.

Let $z$ be the commutator $[v,e]:=ve-ev$.  It is a 1-cycle, whose
homology class~$[z]\in H_1(C)=H_1(F)$ we denote by~$\zeta$.  Let
$\varepsilon=[v]\in H_0(C)=H_0(F)=\Z$.  We will give two versions of the
homology calculation.  The first arises naturally from the method of
proof, but the second exhibits the decomposition of the homology with
respect to right and left multiplication by~$\varepsilon$.

\begin{theorem}
\label{thr:1}
\hspace*{0pt}
\begin{enumerate}
\item
$H_*(F)$ is generated as a graded ring by $\varepsilon$ and~$\zeta$, which
satisfy the relations
\begin{gather*}
\varepsilon^2=\varepsilon\\
\varepsilon\zeta = \zeta - \zeta\varepsilon\;.
\end{gather*}
$H_n(F)$ is free abelian of rank~2 for all $n\geq 1$, with basis
$\zeta^n,\zeta^n \varepsilon$.
\item
Let $\alpha=\zeta\varepsilon$ and $\beta=\varepsilon\zeta$.  Then $H_*(F)$ is
generated as a ring by $\varepsilon$, $\alpha$, and~$\beta$, which satisfy
the relations
\begin{gather*}
\varepsilon^2=\varepsilon\\
\varepsilon\alpha=\beta \varepsilon=0\\
\alpha \varepsilon =\alpha\;, \quad \varepsilon\beta=\beta\;.
\end{gather*}
We have $\alpha^2=\beta^2=0$, and the alternating products
$\alpha\beta\alpha\cdots$ and $\beta\alpha\beta\cdots$ form a basis
for~$H_*(F)$ in positive dimensions.
\end{enumerate}
\end{theorem}

\begin{proof}
It will be convenient to adjoin an identity to~$C$, getting a differential
graded ring $R:=\Z\oplus C$ with identity.  It is the ring of all
noncommuting polynomials in~$v,e$ or, equivalently, the tensor algebra
over~\Z\ of the free \Z-module generated by $v$ and~$e$.  We will
compute~$H_*(R)$, which is simply~$H_*(F)$ with an extra summand~\Z\
in dimension~0.

Let $w=v^2-v$ and let $S<R$ be the subring (with identity) generated
by $e$, $w$, and~$z$.  I claim that $S$ is the ring of noncommuting
polynomials in the variables $e,w,z$ and that $R=S\oplus Sv$.  In
other words, the set~\M\ consisting of the monomials in~$e,w,z$ and
their products with~$v$ forms a basis for~$R$ as a \Z-module.  We will
show that \M\ spans~$R$; a counting argument that can be found in
\cite{stein92:_group}*{p.~498} then shows that \M\ is a basis.  It
suffices to show that any monomial in~$v,e$ can be rewritten as a
linear combination of the monomials in~\M.  We may assume that $v$
occurs in the given monomial and that the first occurence of~$v$ is
not at the end.  It is therefore followed either by another~$v$ or
by~$e$.  If it is followed by~$v$, then we replace the resulting~$v^2$
by $w+v$; otherwise we replace~$ve$ by $ev+z$.  We now have two terms,
and we repeat the process with each of those.  Continuing in this way,
we arrive after finitely steps at a linear combination of elements
of~\M.  (Termination of the process follows from the fact that we are
moving $v$'s to the right and, when we rewrite~$v^2$, reducing the
total number of~$v$'s.)

We now know that $S$ is the tensor algebra $\Z\oplus D \oplus (D\otimes
D) \oplus\cdots$, where $D$ is the \Z-module generated by $e,w,z$.
Each summand is a chain subcomplex, whose homology can be computed by
the K\"unneth formula.  By inspection, $H_*(D)=\Z$, concentrated in
dimension~1 and generated by~$\zeta$; so we conclude that $H_n(S)=\Z$
for every $n\geq 0$, generated by~$\zeta^n$.  The complementary
summand~$Sv$ is isomorphic to~$S$ as a chain complex, so it
contributes a second~\Z, generated by~$\zeta^n\varepsilon$, for each~$n$.
To complete the proof of~(1), we need to check the relations.  The
relation $\varepsilon^2=\varepsilon$ is immediate from $\partial e = v^2 - v$.
For the second, one checks that $\partial (e^2)$ is the commutator
$[v^2-v,e]=[v^2,e]-z=vz+zv-z$, whence $\varepsilon\zeta+\zeta\varepsilon-\zeta=0$.

Turning now to (2), the stated relations are immediate from the
relations in~(1).  For example, we get $\varepsilon\zeta\varepsilon=0$, which says
precisely that $\varepsilon\alpha=\beta\varepsilon=0$, by multiplying the relation
$\varepsilon\zeta = \zeta - \zeta\varepsilon$ by~$\varepsilon$ on the right.  The
equation $\alpha^2=0$ now follows by computing $\alpha\varepsilon\alpha$ in
two different ways, and similarly for $\beta^2=0$.  Since
$\zeta=\alpha+\beta$, it is clear that $H_*(R)$ is generated by
$\varepsilon,\alpha,\beta$ and hence that the alternating products of
$\alpha$ and~$\beta$ generate~$H_n(R)$ additively for $n>0$.  Since we
already know that $H_n(R)\iso \Z\oplus\Z$, we conclude that the
alternating products form a basis.  (Alternatively, note that
$\zeta^n=(\alpha+\beta)^n$ is the sum of the two alternating products
of length~$n$, while $\zeta^n\varepsilon$ is the one that ends in~$\alpha$.)
\end{proof}

\begin{remark}
\label{rem:3}
In both (1) and (2), it is easy to show that the stated relations form
a system of defining relations.
\end{remark}
  
\begin{remark}
\label{rem:4}
The subring $\varepsilon H_*(F) \varepsilon < H_*(F)$, on which
$\varepsilon$ is a 2-sided identity, is easily checked to be the
polynomial ring~$\Z[t]$ generated by $t:=\beta\alpha\in H_2(F)$.  I
claim that it can be identified with~$H_*(F')$, where $F'$ is the
commutator subgroup of~$F$.  To deduce this from what we have already
done, recall first that $F'$ is the subgroup of~$F$ consisting of
dyadic PL-homeomorphisms with support in the interior of~$[0,1]$ (see
Section~\ref{sec:abelianization-f}).  From this one sees that $F'$ is
closed under the product~$\mu$ defined in
Section~\ref{sec:definition}, so $H_*(F')$ is a ring.  Moreover, it is
not hard to see that the canonical generator of $H_0(F')$ is a 2-sided
identity, so the map $H_*(F')\to H_*(F)$ induced by the inclusion
$F'\into F$ is a ring homomorphism with image contained in~$\Z[t]$.
The next observation is that $F$ behaves homologically as though it
were $F'\times F_{\ab}$; more precisely, there is a homomorphism $F\to
F'\times F_{\ab}$ that induces an isomorphism in homology (see
\cite{stein92:_group}*{Lemma~4.1}).  One can deduce that $H_*(F')$
maps injectively onto a direct summand of $H_*(F)$, and a counting
argument [based on $H_*(F)\iso H_*(F')\otimes H_*(F_{\ab})$] shows
that the image is infinite cyclic in every even dimension.  Since the
image is a direct summand of~$\Z[t]$, it must equal~$\Z[t]$.
\end{remark}

\section{Calculation of the cup product}
\label{sec:calc-cup-prod}

By standard arguments, as in the homology theory of topological
groups, the diagonal map $F\to F\times F$ (or $X\to X\times X$)
induces a ring homomorphism
\[
\Delta\colon H_*(F) \to H_*(F) \otimes H_*(F)\;,
\]
whose dual is the cup product in cohomology.  Thus the cup product
will be known if we can compute $\Delta$ on a set of ring generators
of~$H_*(F)$, such as $\varepsilon,\alpha,\beta$.  By general
principles, we know that $\Delta(\varepsilon) = \varepsilon \otimes
\varepsilon$ and that $\alpha$ and~$\beta$ are primitive, i.e.,
$\Delta(\alpha) = \alpha\otimes\varepsilon + \varepsilon\otimes\alpha$
and similarly for~$\beta$.  One can deduce, after some calculations,
that the cohomology ring is given by formula~\eqref{eq:1} in the
introduction.  The calculations are simpler if we make use of the
homology equivalence $F\to F'\times F_{\ab}$ mentioned in
Remark~\ref{rem:4} above, so we will phrase the statement and proof in
those terms.

Recall that the divided polynomial ring $\Gamma(u)$ on one
generator~$u$ is the subring of the polynomial ring~$\Q[u]$ generated
additively by the elements $u^{(i)} := u^i/i!$
\thinspace ($i\geq 0$).  Note that
\begin{equation}
\label{eq:4}
u^{(i)} u^{(j)} = \binom{i+j}{i} u^{(i+j)}\;,
\end{equation}
so the \Z-span of the $u^{(i)}$ is indeed a ring.

\begin{theorem}
\label{thr:2}
\hspace*{0pt}
\begin{enumerate}
\item
$H^*(F') \iso \Gamma(u)$, where $\deg u = 2$.
\item
There are ring isomorphisms
\[
H^*(F) \iso H^*(F')\otimes H^*(F_{ab}) \iso \Gamma(u) \otimes
\bigwedge(a,b)\;,
\]
where $\deg a = \deg b = 1$ and $\deg u = 2$.
\end{enumerate}
\end{theorem}

\begin{proof}
Recall that $H_*(F') = \Z[t]$.  By general principles (or a
calculation of $\Delta(t)=\Delta(\beta)\Delta(\alpha)$), the
element~$t\in H_2(F')$ is primitive, so
\[
\Delta(t^n)=(t \otimes 1 + 1\otimes t)^n = \sum_{i+j=n} \binom{i+j}{i}
t^i \otimes t^j\;.
\]
Dualizing, we obtain a basis $(u^{(n)})_{n\geq 0}$ for $H^*(F')$, with
multiplication law as in~\eqref{eq:4}.  This proves~(1), and (2)
follows at once.
\end{proof}


\begin{bibdiv}
\begin{biblist}

\input{xrefs.ltb}
\bib{belk04:_f}{thesis}{
  author={Belk, James M.},
  title={Thompson's group $F$},
  type={Ph.D. Thesis},
  date={2004},
  address={Cornell University},
}

\bib{belk:_fores_f}{article}{
  author={Belk, James~M.},
  author={Brown, Kenneth~S.},
  title={{Forest diagrams for elements of Thompson's group~$F$}},
  date={2004},
  journal={Internat. J. Algebra Comput.},
  status={to appear},
}

\bib{brown87:_finit}{article}{
  author={Brown, Kenneth S.},
  title={Finiteness properties of groups},
  booktitle={Proceedings of the Northwestern conference on cohomology of groups (Evanston, Ill., 1985)},
  journal={J. Pure Appl. Algebra},
  volume={44},
  date={1987},
  number={1-3},
  pages={45\ndash 75},
  issn={0022-4049},
  review={\MR {MR885095 (88m:20110)}},
}

\bib{brown92}{article}{
  author={Brown, Kenneth S.},
  title={The geometry of finitely presented infinite simple groups},
  book={MR1230626},
  pages={121\ndash 136},
  review={\MR {MR1230631 (94f:20059)}},
}

\bib{brown94:_cohom}{book}{
  author={Brown, Kenneth S.},
  title={Cohomology of groups},
  series={Graduate Texts in Mathematics},
  volume={87},
  note={Corrected reprint of the 1982 original},
  publisher={Springer-Verlag},
  place={New York},
  date={1994},
  pages={x+306},
  isbn={0-387-90688-6},
  review={\MR {MR1324339 (96a:20072)}},
}

\bib{brown83:_fp_hnn}{article}{
  author={Brown, Kenneth S.},
  author={Geoghegan, Ross},
  title={${\rm FP}\sb {\infty }$ groups and HNN extensions},
  journal={Bull. Amer. Math. Soc. (N.S.)},
  volume={9},
  date={1983},
  number={2},
  pages={227\ndash 229},
  issn={0273-0979},
  review={\MR {MR707963 (85a:20024)}},
}

\bib{brown84:_fp}{article}{
  author={Brown, Kenneth S.},
  author={Geoghegan, Ross},
  title={An infinite-dimensional torsion-free ${\rm FP}\sb {\infty }$ group},
  journal={Invent. Math.},
  volume={77},
  date={1984},
  number={2},
  pages={367\ndash 381},
  issn={0020-9910},
  review={\MR {MR752825 (85m:20073)}},
}

\bib{brown85:_cohom}{article}{
  author={Brown, Kenneth S.},
  author={Geoghegan, Ross},
  title={Cohomology with free coefficients of the fundamental group of a graph of groups},
  journal={Comment. Math. Helv.},
  volume={60},
  date={1985},
  number={1},
  pages={31\ndash 45},
  issn={0010-2571},
  review={\MR {MR787660 (87b:20066)}},
}

\bib{cannon96:_introd_richar}{article}{
  author={Cannon, James W.},
  author={Floyd, William J.},
  author={Parry, Walter R.},
  title={Introductory notes on Richard Thompson's groups},
  journal={Enseign. Math. (2)},
  volume={42},
  date={1996},
  number={3-4},
  pages={215\ndash 256},
  issn={0013-8584},
  review={\MR {MR1426438 (98g:20058)}},
}

\bib{clark69:_variet}{article}{
  author={Clark, David M.},
  title={Varieties with isomorphic free algebras},
  journal={Colloq. Math.},
  volume={20},
  date={1969},
  pages={181\ndash 187},
  review={\MR {MR0244128 (39 \#5445)}},
}

\bib{dydak77:_fanr}{article}{
  author={Dydak, Jerzy},
  title={A simple proof that pointed FANR-spaces are regular fundamental retracts of ANR's},
  language={English, with Russian summary},
  journal={Bull. Acad. Polon. Sci. S\'er. Sci. Math. Astronom. Phys.},
  volume={25},
  date={1977},
  number={1},
  pages={55\ndash 62},
  review={\MR {MR0442918 (56 \#1293)}},
}

\bib{dydak77}{article}{
  author={Dydak, Jerzy},
  title={1-movable continua need not be pointed 1-movable},
  language={English, with Russian summary},
  journal={Bull. Acad. Polon. Sci. S\'er. Sci. Math. Astronom. Phys.},
  volume={25},
  date={1977},
  number={6},
  pages={559\ndash 562},
  review={\MR {MR0474226 (57 \#13873)}},
}

\bib{freyd93:_split}{article}{
  author={Freyd, Peter},
  author={Heller, Alex},
  title={Splitting homotopy idempotents. II},
  journal={J. Pure Appl. Algebra},
  volume={89},
  date={1993},
  number={1-2},
  pages={93\ndash 106},
  issn={0022-4049},
  review={\MR {MR1239554 (95h:55015)}},
}

\bib{higman74:_finit}{book}{
  author={Higman, Graham},
  title={Finitely presented infinite simple groups},
  publisher={Department of Pure Mathematics, Department of Mathematics, I.A.S. Australian National University, Canberra},
  date={1974},
  pages={vii+82},
  review={\MR {MR0376874 (51 \#13049)}},
}

\bib{jonsson61}{article}{
  author={J{\'o}nsson, Bjarni},
  author={Tarski, Alfred},
  title={On two properties of free algebras},
  journal={Math. Scand.},
  volume={9},
  date={1961},
  pages={95\ndash 101},
  review={\MR {MR0126399 (23 \#A3695)}},
}

\bib{quillen73:_higher_k}{article}{
  author={Quillen, Daniel},
  title={Higher algebraic $K$-theory. I},
  book={MR0325307},
  pages={85\ndash 147},
  review={\MR {MR0338129 (49 \#2895)}},
}

\bib{quillen78:_homot}{article}{
  author={Quillen, Daniel},
  title={Homotopy properties of the poset of nontrivial $p$-subgroups of a group},
  journal={Adv. in Math.},
  volume={28},
  date={1978},
  number={2},
  pages={101\ndash 128},
  issn={0001-8708},
  review={\MR {MR493916 (80k:20049)}},
}

\bib{smirnov71:_cantor}{article}{
  author={Smirnov, D. M.},
  title={Cantor algebras with one generator. I},
  language={Russian},
  journal={Algebra i Logika},
  volume={10},
  date={1971},
  pages={61\ndash 75},
  review={\MR {MR0296006 (45 \#5067)}},
}

\bib{stein92:_group}{article}{
  author={Stein, Melanie},
  title={Groups of piecewise linear homeomorphisms},
  journal={Trans. Amer. Math. Soc.},
  volume={332},
  date={1992},
  number={2},
  pages={477\ndash 514},
  issn={0002-9947},
  review={\MR {MR1094555 (92k:20075)}},
}

\bib{swierczkowski61}{article}{
  author={{\'S}wierczkowski, Stanis{\l }aw},
  title={On isomorphic free algebras},
  journal={Fund. Math.},
  volume={50},
  date={1961/1962},
  pages={35\ndash 44},
  review={\MR {MR0138573 (25 \#2017)}},
}




\end{biblist}
\end{bibdiv}
\end{document}